\documentclass{article}
\usepackage{amsmath,amsthm,amsfonts,amssymb, latexsym}

\newtheorem{thm}{Theorem}

\newtheorem{lem}[thm]{Lemma}
\newtheorem{prop}[thm]{Proposition}
\newtheorem{cor}[thm]{Corollary}

\theoremstyle{definition}
\newtheorem*{defin}{Definition}

\DeclareMathOperator{\vol}{Vol} \DeclareMathOperator{\re}{Re}
\DeclareMathOperator{\im}{Im} \DeclareMathOperator{\su}{SU}
\DeclareMathOperator{\gl}{GL} \DeclareMathOperator{\symp}{Sp}
\DeclareMathOperator{\un}{U} \DeclareMathOperator{\so}{SO}
\DeclareMathOperator{\slin}{SL}
\newcommand{\dsum}{\displaystyle \sum}

\def\rr{\mathbb{R}}
\def\cc{\mathbb{C}}
\def\zz{\mathbb{Z}}

\def\too{\longrightarrow}

\def\M{M_\alpha}
\def\DD{\mathrm{D}}
\def\dd{\mathrm{d}}

\begin{document}

\title{Connected Sums of Special Lagrangian Submanifolds}
\author{Dan A.\ Lee\\ Stanford University\\ dalee@math.stanford.edu}
\date{\today}
\maketitle

\begin{abstract}
Let $M_1$ and $M_2$ be special Lagrangian submanifolds of a
compact Calabi-Yau manifold $X$ that intersect transversely at a
single point. We can then think of $M_1\cup M_2$ as a singular
special Lagrangian submanifold of $X$ with a single isolated
singularity. We investigate when we can regularize $M_1\cup M_2$
in the following sense: There exists a family of Calabi-Yau
structures $X_\alpha$ on $X$ and a family of special Lagrangian
submanifolds $M_\alpha$ of $X_\alpha$ such that $M_\alpha$
converges to $M_1\cup M_2$ and $X_\alpha$ converges to the
original Calabi-Yau structure on $X$.  We prove that a
regularization exists in two key cases: (1) when $\dim_\cc X=3$,
Hol$(X)=\su(3)$, and $[M_1]$ is not a multiple of $[M_2]$ in
$H_3(X)$, and (2) when $X$ is a torus with $\dim_{\cc}X\geq3$,
$M_1$ is flat, and the intersection of $M_1$ and $M_2$ satisfies a
certain angle criterion. One can easily construct examples of the
second case, and thus as a corollary we construct new examples of
non-flat special Lagrangian submanifolds of Calabi-Yau tori.
\end{abstract}

\section{Introduction}

\indent One of the fundamental problems in special Lagrangian
geometry is to understand moduli spaces of special Lagrangian
submanifolds (SLags).  Much interest in this problem arises from
the study of mirror symmetry since it is related to the SYZ
Conjecture \cite{stryauzas}.  McLean's deformation theorem
\cite{mclean} together with some work by Hitchin \cite{hitchin}
provide some understanding of these moduli spaces locally near
nonsingular SLags, but in order to understand these moduli spaces
globally, we need to understand singular SLags.  Recently, some
research has focused on the more modest goal of understanding
SLags with isolated conical singularities.  For example, see
\cite{joyceconical1,joyceconical2,joyceconical3,joyceconical4}. In
order to study SLags with isolated conical singularities, we need
to know something about the SLag cones in $\cc^n$ on which these
singularities are modelled.  Some researchers have undertaken the
study of these SLag cones.  See \cite{haskins}.

In this paper we restrict our attention to isolated conical
singularities modelled on a very simple type of SLag cone, namely,
the union of two transversely intersecting SLag planes in $\cc^n$.
By the work of Lawlor \cite{lawlor}, we know that such a SLag cone
can be deformed through a family of nonsingular SLags in $\cc^n$
if the two planes meet a certain angle criterion to be described
later. (This criterion is always satisfied when $n\leq3$.)  This
local regularization holds out hope that if we have a singular
SLag with this simple type of singularity, then it can be globally
regularized, that is, it can be deformed through a family of
nonsingular SLags.  This simple type of singularity arises when a
connected, immersed SLag intersects itself, and when two embedded
SLags intersect.  The first case has already been treated by
Yng-Ing Lee \cite{ylee} who answered the question in the
affirmative:  A compact, connected, immersed SLag with an isolated
point of transverse self-intersection satisfying the angle
criterion can be regularized.  The second case is more difficult,
and that is the case which we consider in this paper. Simply put,
our problem is to try to regularize the union of two compact
embedded SLags with an isolated point of transverse intersection
satisfying the angle criterion.  A problem related to ours has
been solved by Butscher \cite{butscher}:  The union of two
embedded SLags with boundary in $\cc^n$ with $n\geq3$ with an
isolated point of transverse intersection satisfying the angle
criterion can be regularized.  In Butscher's paper, the
regularization takes advantage of the freedom to deform the
boundary of the singular SLag.  In our problem, we have no
boundaries, and therefore we cannot use the added degrees of
freedom.  Indeed, our problem as stated above probably cannot be
solved.  We must introduce another degree of freedom, and we do
this by deforming the Calabi-Yau structure of the ambient
manifold.

Before we state our results, we recall some basic definitions and
facts.

\begin{defin}
A \textbf{Calabi-Yau structure} (or CY structure) on a compact
$2n$-fold $X$ is a 3-tuple $(J,\omega,\Omega)$ such that $J$ is a
complex structure on $X$, $\omega$ is a K\"{a}hler form with
respect to $J$, and $\Omega$ is a holomorphic $(n,0)$ form with
respect to $J$ such that
\begin{equation}\label{normalization}
{\omega^n\over
n!}=(-1)^{{1\over2}n(n-1)}\left({i\over2}\right)^n\Omega\wedge\overline{\Omega}
\end{equation}
It is a fact that $\re\Omega$ is a calibration with respect to the
K\"{a}hler metric.  We say that a submanifold $M$ of $X$ is
\textbf{special Lagrangian} iff $M$ is calibrated by $\re\Omega$.
\end{defin}

The special Lagrangian condition on $M$ is equivalent to the
vanishing of both $\omega$ and $\im\Omega$ on $M$.  If $(X,J)$
admits a Calabi-Yau structure at all, then in each K\"{a}hler
class there is a unique K\"{a}hler form $\omega$ such that
$(X,J,\omega)$ admits a Calabi-Yau structure.  (In this case, the
K\"{a}hler metric corresponding to $\omega$ is Ricci-flat.)  Also,
if $(X,J)$ admits a Calabi-Yau structure, then $\Omega$ is
uniquely determined up to a complex constant.  If we also choose
$\omega$, then the normalization (\ref{normalization}) uniquely
determines $\Omega$ up to a phase. Because of these facts, a
choice of Calabi-Yau structure amounts to a choice of complex
structure, a K\"{a}hler class, and a phase. For more general
background on Calabi-Yau manifolds and special Lagrangian
geometry, see \cite{schoensurvey,joycesurvey,joycebook,harvey}.

We now describe the angle criterion.
\begin{defin}
Given any two $n$-dimensional oriented linear subspaces $\eta$ and
$\xi$ of $\rr^{2n}$, there exist \textbf{characterizing angles}
$\theta_1,\ldots,\theta_n$, together with an orthonormal basis
$e_1,\ldots,e_{2n}$ of $\rr^{2n}$, such that
$$0\leq\theta_1\leq\cdots\leq\theta_{n-1}\leq{\pi\over 2}\text{ and }
\theta_{n-1}\leq\theta_n\leq\pi-\theta_{n-1}$$ while
$$\eta=e_1\wedge\ldots\wedge e_n$$
and
$$\xi=[(\cos\theta_1)e_1+(\sin\theta_1)e_{n+1}]\wedge\ldots\wedge
[(\cos\theta_n)e_n+(\sin\theta_n)e_{2n}].$$ We say that $\eta$ and
$\xi$ satisfy the \textbf{angle criterion} iff the characterizing
angles between $\eta$ and $-\xi$ satisfy $\sum_{i=1}^n \theta_i
=\pi$.
\end{defin}
The Lawlor-Nance Angle Theorem states that a pair of oriented
planes $(\eta,\xi)$ is minimizing iff the characterizing angles
between $\eta$ and $-\xi$ satisfy \mbox{$\sum_{i=1}^n \theta_i
\geq\pi$}. Therefore the angle criterion may be thought of as
describing the ``borderline case" of minimizing pairs of planes.
See \cite{harvey} for more on characterizing angles and the Angle
Theorem.

We are now ready to state our main theorem.
\begin{thm}[Main Theorem]
Let $M_1$ and $M_2$ be two embedded special Lagrangian
submanifolds of a Calabi-Yau manifold $(X,J,\omega,\Omega)$ such
that $n=\dim_\cc X\geq3$ and the holonomy of the K\"{a}hler metric
is exactly $\su(n)$.  Assume that $M_1$ and $M_2$ intersect
transversely at a single point $p$ such that the tangent cone of
$M_1\cup M_2$ at $p$ satisfies the angle criterion.  Further
assume that $\re[H^{n-1,1}(X)\oplus H^{1,n-1}(X)]$ is not
contained in the kernel of
${[M_1]\over\vol(M_1)}-{[M_2]\over\vol(M_2)}$, thought of as a
functional on $H^n(X)$.  Then there exists a family of Calabi-Yau
structures $(J_\alpha,\omega_\alpha, \Omega_\alpha)$ on $X$
converging to $(J,\omega,\Omega)$ and a family of embedded
submanifolds $M_\alpha\subset X$ converging to $M_1\cup M_2$ such
that $M_\alpha$ is special Lagrangian in
$(X,J_\alpha,\omega_\alpha, \Omega_\alpha)$.
\end{thm}

When $n=3$, the theorem reduces to the following nice result.
\begin{cor}
Let $M_1$ and $M_2$ be two embedded special Lagrangian
submanifolds of a Calabi-Yau manifold $(X,J,\omega,\Omega)$ such
that $\dim_\cc X=3$ and the holonomy of the K\"{a}hler metric is
exactly $\su(3)$.  Assume that $M_1$ and $M_2$ intersect
transversely at a single point and that $[M_1]$ is not a multiple
of $[M_2]$ in $H_n(X)$. Then there exists a family of Calabi-Yau
structures $(J_\alpha,\omega_\alpha, \Omega_\alpha)$ on $X$
converging to $(J,\omega,\Omega)$ and a family of embedded
submanifolds $M_\alpha\subset X$ converging to $M_1\cup M_2$ such
that $M_\alpha$ is special Lagrangian in
$(X,J_\alpha,\omega_\alpha, \Omega_\alpha)$.
\end{cor}
\begin{proof}
As mentioned earlier, the angle criterion is automatically
satisfied when $n=3$.  It suffices to show that the homology
condition in the Corollary implies the one in the Main Theorem.
Since $[\Omega]$ spans $H^{3,0}(X)$, it follows that $[\re\Omega]$
and $[\im\Omega]$ span $H^{3,0}(X)\oplus H^{0,3}(X)$.  By
assumption, ${[M_1]\over\vol(M_1)}-{[M_2]\over\vol(M_2)}\neq0$,
and by the special Lagrangian condition, this real homology class
is zero on both $[\re\Omega]$ and $[\im\Omega]$.  So there exists
some element of $\re[H^{2,1}(X)\oplus H^{1,2}(X)]$ which is not in
the kernel of ${[M_1]\over\vol(M_1)}-{[M_2]\over\vol(M_2)}$.
\end{proof}

Since a CY torus has trivial holonomy, our Main Theorem does not
apply to this important case. However, we can still prove a
version of the theorem in this setting.
\begin{thm}[Torus Version]
Let $M_1$ and $M_2$ be two embedded special Lagrangian
submanifolds of a Calabi-Yau torus $(T,J,\omega,\Omega)$ such that
$\dim_\cc T\geq 3$ and $M_1$ is flat.  Assume that $M_1$ and $M_2$
intersect transversely at a single point $p$ such that the tangent
cone of $M_1\cup M_2$ at $p$ satisfies the angle criterion. Then
there exists a family of Calabi-Yau structures $(J_\alpha,\omega,
\Omega_\alpha)$ on $T$ converging to $(J,\omega,\Omega)$ and a
family of embedded submanifolds $M_\alpha\subset T$ converging to
$M_1\cup M_2$ such that $M_\alpha$ is special Lagrangian in
$(T,J_\alpha,\omega, \Omega_\alpha)$.
\end{thm}
The Main Theorem and the Torus Version share the hypothesis that
$M_1$ and $M_2$ must intersect at a single point, but this
condition is somewhat artificial. In light of the proof to follow,
as long as there exists an \emph{isolated} transverse intersection
point $p$, we can still regularize the singularity at $p$, but the
$\M$'s will only be immersed rather than embedded. However, if
$M_1\cap M_2$ is a finite set of isolated transverse intersection
points, all of which satisfy the angle criterion, then we can
recover the embeddedness as follows: We first apply our result to
one of these intersection points,\footnote{If either $M_1$ or
$M_2$ is not connected, then we apply our result multiple times.}
and then we apply Yng-Ing Lee's result on immersed SLags to each
of the other intersection points. This procedure is possible
because the property of being a transverse intersection point
satisfying the angle criterion is an open condition with respect
to the relevant topology.  See Lemma \ref{anglecrit}.

It is a simple matter to construct infinitely many distinct pairs
of flat SLag tori satisfying the angle criterion in the standard
CY torus, $\cc^n/\zz^{2n}$. Applying the discussion in the
previous paragraph, we immediately obtain the following result.
\begin{cor}
There exist non-flat embedded special Lagrangian submanifolds of
Calabi-Yau tori.
\end{cor}

With some extra work, the methods of this paper can probably be
used to prove that our results hold in dimension two also.  The
methods of algebraic geometry should also apply in the dimension
two case.

Concurrent with the writing of this paper, Joyce has produced some
results on the general problem of desingularizing special
Lagrangians with isolated conical singularities in almost
Calabi-Yau manifolds \cite{joyceconical3,joyceconical4}. In
particular, Theorem 7.11 of \cite{joyceconical4} combined with
Lemma \ref{deform} of this paper and an understanding of the
Lawlor necks can be used to prove the results of this paper.  Note
that Lemma \ref{deform} is the main ingredient of the Key Lemma of
this paper. The methods used by Joyce are different from those
presented here, and because of the added generality, the proofs
are also more complicated.
\\
\\
\noindent {\sc{Acknowledgements:}} I would like to thank Rick
Schoen for suggesting the problem, listening to my ideas, and
offering many helpful suggestions.  This research was partially
supported by a NSF Graduate Research Fellowship.

\section{Preliminaries}

\indent  The Main Theorem and the Torus Version share certain
assumptions: We have two embedded special Lagrangian submanifolds
 $M_1$ and $M_2$ of a Calabi-Yau manifold $(X,J,\omega,\Omega)$
 with $\dim_\cc X\geq 3$.  We also assume that $M_1$
and $M_2$ intersect transversely at a single point $p$ such that
the tangent cone of $M_1\cup M_2$ at $p$ satisfies the angle
criterion.  This is the situation we assume from now until the
proofs of the Key Lemma, which will depend on the additional
assumptions in the two cases.  We also assume without loss of
generality that $M_1$ and $M_2$ are connected.

We now explain the idea behind these results. We wish to construct
a family of approximate solutions $\M$ such that $\M$ converges to
$M=M_1\cup M_2$, $\M$ is exactly Lagrangian, and $\M$ is very
close to being special Lagrangian. Once we have these $\M$'s, we
can construct small Hamiltonian deformations of them and hope that
at least one of them is exactly special Lagrangian. This is
actually too much to hope for, but we can add another degree of
freedom to this deformation by simultaneously deforming the
Calabi-Yau structure $(J_t,\omega_t,\Omega_t)$ and $\M$ itself so
that the deformations of $\M$ remain Lagrangian with respect to
$\omega_t$. Using these deformations we define a deformation
operator whose solutions correspond to special Lagrangians in
$(X,J_t,\omega_t,\Omega_t)$. Using the Inverse Function Theorem
together with some estimates, we obtain the desired solutions. The
work lies in obtaining the appropriate estimates.

First we construct our family of approximate solutions $\M$.  This
is where the angle criterion is relevant. Given
\mbox{$\phi_1,\ldots,\phi_n\in\rr$}, we define
$P(\phi_1,\ldots,\phi_n)$ to be the oriented plane
$[(\cos\phi_1){\partial\over\partial x^1}+(\sin\phi_1)
{\partial\over\partial y^1}]\wedge\ldots\wedge
[(\cos\phi_n){\partial\over\partial x^n}+(\sin\phi_n)
{\partial\over\partial y^n}]$.  By the work of Lawlor
\cite{lawlor}, we know that for any
$\theta_1,\ldots,\theta_n\in(0,\pi)$ satisfying
$\sum_{j=1}^n\theta_j=\pi$, there exists a special Lagrangian
submanifold $N$ of $\cc^n$ that is asymptotic in an oriented sense
to the two planes $P(0,\ldots,0)$ and
$-P(-\theta_1,\ldots,-\theta_n)$. This submanifold $N$ has the
property that $\epsilon N$ converges to $[P(0,\ldots,0)]\cup
[-P(-\theta_1,\ldots,-\theta_n)]$ in an appropriate sense as
$\epsilon\to 0$. These $N$'s, as well as their images under
$\su(n)\times$(dilations) are called \textbf{Lawlor necks}. Note
that there is also a Lawlor neck asymptotic in an oriented sense
to the two planes $P(0,\ldots,0)$ and
$-P(\theta_1,\ldots,\theta_n)$.
\begin{lem}\label{anglecrit}
If $\eta$ and $\xi$ are two special Lagrangian planes in $\cc^n$,
then there exists a Lawlor neck asymptotic in an oriented sense to
$\eta$ and $\xi$ if and only if $\eta$ and $\xi$ are transverse
planes satisfying the angle criterion. Moreover, both of these
equivalent conditions are open conditions in the space of pairs of
special Lagrangian planes. Finally, when $n\leq3$, every pair of
special Lagrangian planes satisfies the angle criterion.
\end{lem}
\begin{proof}
Let $\eta$ and $\xi$ be two transverse special Lagrangian planes
in $\cc^n$.  The necessity of the angle criterion for the
existence of a Lawlor neck is trivial, so we need only prove
sufficiency.  Assume that $\eta$ and $\xi$ satisfy the angle
criterion.
 Without loss of generality, we may assume that
$\eta=P(0,\ldots,0)$ by performing an $\su(n)$ change of
coordinates.  Since $-\xi$ is Lagrangian and transverse to $\eta$,
we have $-\xi=P(\phi_1,\ldots,\phi_n)$ for some
$\phi_1,\ldots,\phi_n\in(-\pi,0)\cup(0,\pi)$, after an $\so(n)$
change of coordinates, where $\so(n)\subset\su(n)$ is the standard
inclusion. See \cite{butscherthesis} for a proof of this fact.
Note that we still have $\eta=P(0,\ldots,0)$ in this coordinate
system.  However, the unordered list of angles,
$\phi_1,\ldots,\phi_n$, is not uniquely determined. Adding or
subtracting $\pi$ to any of the $\phi_j$'s merely changes the
orientation, therefore doing this an even number of times leaves
the oriented plane unchanged. We can find a canonical description
by placing more restrictions on the angles. We demand that at most
one of the $\phi_j$'s has $|\phi_j|>{\pi\over 2}$, and if there is
such a $|\phi_{j_0}|$, we demand that $|\phi_{j_0}|\leq
\pi-|\phi_j|$ for all $j$.  It is easy to verify that this
prescription gives us a new unordered list of angles. Now observe
that $|\phi_1|,\ldots,|\phi_n|$ is precisely the unordered list of
characterizing angles for the planes $\eta$ and $-\xi$, therefore
$\sum_{j=1}^n |\phi_j|=\pi$ by assumption. Since $\xi$ is special
Lagrangian, $\sum_{j=1}^n\phi_j\equiv\pi$ (mod $2\pi$). So we have
$\sum_{j=1}^n\phi_j=\pm\pi$.  In the positive case, we have
$\sum_{j=1}^n (|\phi_j|-\phi_j)=0$. Therefore all of the
$\phi_j$'s are positive and their sum is $\pi$. Similarly, in the
negative case, we see that all of the $\phi_j$'s are negative and
their sum is $-\pi$.  In either case, our brief discussion of
Lawlor necks above shows that there exists a Lawlor neck
asymptotic in an oriented sense to $\eta$ and $-\xi$.

We now turn to the second statement in the Lemma.  Let $\eta$ and
$\xi$ be two transverse special Lagrangian planes satisfying the
angle criterion.  By the previous discussion, there exist $\su(n)$
coordinates in which $\eta=P(0,\ldots,0)$ and
$-\xi=P(\theta_1,\ldots,\theta_n)$, where
$\theta_1,\ldots,\theta_n\in(0,\pi)$ and $\sum_{j=1}^n
\theta_j=\pi$. (Or we have the negative case which is similar.)
The key observation is that if $\xi'$ is a SLag plane close to
$\xi$, then we can make an $\so(n)\subset\su(n)$ change of
coordinates so that $-\xi'=P(\theta'_1,\ldots,\theta'_n)$, where
$\theta'_1,\ldots,\theta'_n$ is close to
$\theta_1,\ldots,\theta_n$ in the natural topology on unordered
lists of $n$ objects. This essentially follows from the fact that
the unordered list of eigenvalues (with algebraic multiplicity) of
a matrix continuously depends on the matrix.  Since the change of
coordinates was $\so(n)$, we still have $\eta=P(0,\ldots,0)$ in
this coordinate system.  Since $\xi'$ is special Lagrangian,
$\sum_{j=1}^n \theta'_j\equiv\pi$ (mod $2\pi$), so for $\xi'$
close enough to $\xi$, we see that $\sum_{j=1}^n \theta'_j=\pi$
and $\theta'_1,\ldots,\theta'_n\in(0,\pi)$, thus proving the
second statement of the Lemma.

The final statement in the Lemma about $n\leq3$ is simple to
verify.
\end{proof}

Now we must use the existence of the local regularization to
produce an approximate global regularization.  The details of this
construction are described in \cite{butscher,butscherthesis,ylee}.
Here we only give a broad overview. Near the singular point $p$,
we can choose a Darboux and normal coordinate system in a ball $B$
around $p$ such that $p=0$ and $\Omega$ approaches $\dd z$
appropriately as we approach $0$. We know $M=M_1\cup M_2$ becomes
close to the tangent cone at $p$ as we approach $0$.  The tangent
cone must be a union of two planes calibrated by $\re \dd z$.  As
long as this pair of planes meets the angle criterion, there
exists a Lawlor neck $N$ asymptotic to those two planes. For
sufficiently small $\alpha>0$ and certain constants $C_\delta$ and
$C_\epsilon$, choose
$$\delta={\alpha\over C_\delta}\text{ and }\epsilon=
{\alpha^{1+1/n}\over C_\epsilon}.$$ See
\cite{butscher,butscherthesis,ylee} for the definitions of
$C_\delta$ and $C_\epsilon$.\footnote{Throughout this paper
$\alpha$ will be the parameter upon which most of our
constructions depend. Because of this, we will explicitly write
out the $\alpha$ dependence of all of our constants, with the only
exceptions being $\delta$ and $\epsilon$. We will use the letter
$C$ without subscript as a generic constant independent of
$\alpha$ whose value may change even in a single chain of
inequalities.
 For consistency we always use $C$ as an
 upper bound.} These constants are chosen so that the following
 construction works and has the stated properties.

We can cut out a small ball $B_\delta(0)$ from $M$, glue a
rescaled Lawlor neck $\epsilon N$ into $B_{\delta\over2}(0)$, and
then interpolate in the annular region.\footnote{From now on we
will write $B_\delta$ for $B_\delta(0)$ where there is no chance
of confusion.} This gives us $\M$ which we can think of as
$M'_1\cup T_1\cup N'\cup T_2\cup M'_2$, where
$M'_i=M_i-B_\delta(0)$, $N'$ is the rescaled Lawlor neck, and the
$T_i$'s are the interpolated regions connecting $M'_i$ to $N'$. It
is evident that $\M$ converges to $M$ by construction, and since
each Lawlor neck has the topology of a cylinder,
$S^{n-1}\times\rr$, $\M$ is topologically the connected sum of
$M_1$ and $M_2$.  We can choose the interpolation so that $\M$ is
exactly Lagrangian. Since $\M$ is Lagrangian, it is a fact that
that at each point of $\M$, $\Omega|_{\M}=e^{i\theta}\vol_{\M}$
for some $\theta$. We call the multi-valued function $\theta$ the
\textbf{Lagrangian angle function}.  This ``function" has the
property that $J\nabla\theta$ is the mean curvature field $H$.  On
a special Lagrangian submanifold, $e^{i\theta}=1$ and $H=0$. One
can show that $\M$ is approximately special Lagrangian in
following sense \cite{butscher,butscherthesis,ylee}.
\begin{lem}\label{theta}
Given the construction above, the Lagrangian angle function and
the mean curvature field on $\M$ satisfy
\begin{eqnarray*}
&|\sin\theta|_0+\alpha^{\beta}[\sin\theta]_\beta+
 \alpha|\nabla\sin\theta|_0&\leq C\alpha\\
&|1-\cos\theta|_0+ \alpha^{\beta}[\cos\theta]_\beta
 +\alpha|\nabla\cos\theta|_0&\leq C\alpha^2\\
&|H|_0+\alpha^{\beta}[H]_\beta+\alpha|\nabla H|_0&\leq C.
\end{eqnarray*}
Moreover, on $\M-B_{\delta}$, $e^{i\theta}=1$ and $H=0$.
\end{lem}

\section{The Deformation Operator}

\indent Now that we have our approximate solutions $\M$, we can
define the relevant deformation operator. Suppose that we have a
smooth deformation $(J_t,\omega_t,\Omega_t)$ of the CY structure
$(J,\omega,\Omega)=(J_0,\omega_0,\Omega_0)$ such that $\omega_t$
is always cohomologous to $\omega$ and $\langle [M_1]+[M_2], [\im
\Omega_t]\rangle=0$. We will choose the appropriate deformation of
CY structure later in this paper; for now assume that we have
already chosen it. By Moser's Theorem, there exists a smooth path
of diffeomorphisms $\Psi_t$ of $X$ such that
\begin{equation}\label{moser}
 \Psi_t^*
\omega_t=\omega. \end{equation}
 By the Lagrangian Neighborhood Theorem, let $U$ be a tubular neighborhood
of $\M$ symplectomorphic to $T^*\M$ so that we have a projection
map $\pi:U\too\M$.  Let $\tau$ be a smooth cutoff function
supported in $U$ such that $\tau = 1$ on ${1\over 2}U$, where
${1\over 2}U$ is defined using the structure of $T^*M$. Observe
that we can choose $U$ to have width greater than ${\epsilon\over
C}$ over $\M\cap B_{\epsilon r_0}$ for some $r_0$ and width
greater than ${1\over C}$ over $M_1'\cup M_2'$, with an inverse
linear interpolation in between.  Now extend any function $h\in
C^{2,\beta}(\M)$ to a function $\tilde{h}\in C^{2,\beta}(X)$ by
defining $\tilde{h}(q)=\tau(q)h(\pi(q))$ on $U$ and $\tilde{h}=0$
outside $U$.  Now define $\Phi_h$ to be the symplectomorphism
generated by the Hamiltonian function $\tilde{h}$.

\begin{defin}
The \textbf{deformation operator}
$F_\alpha:C^{2,\beta}(\M)\times\rr\too C^{0,\beta}(\M)$ is defined
by
$$F_\alpha(h,t)=\langle (\Psi_t\circ\Phi_h)^* (\im\Omega_t),\vol_{\M}
\rangle_{\M}$$ where the metric on $\M$ is the one induced by the
K\"{a}hler metric on $(X,J,\omega)$, independent of $t$.
\end{defin}
Since $\M$ is a Lagrangian submanifold of $(X,\omega)$ and
$\Phi_h$ is a symplectomorphism it follows that $\Phi_h(\M)$ is a
Lagrangian submanifold of $(X,\omega)$. Then by (\ref{moser}), it
follows that $(\Psi_t\circ \Phi_h)(\M)$ is a Lagrangian
submanifold of $(X,\omega_t)$. Clearly, $F_\alpha(h,t)=0$ iff
$\im\Omega_t$ restricted to $(\Psi_t\circ\Phi_h)(\M)$ is
identically zero. Therefore a solution of the equation
$F_\alpha(h,t)=0$ corresponds to a special Lagrangian submanifold
of $(X,J_t,\omega_t, \Omega_t)$, and a small solution corresponds
to a nearby special Lagrangian.  So our goal is to show that for
sufficiently small $\alpha$, $F_\alpha$ has a small solution.  Our
method of constructing such solutions is the following version of
the Inverse Function Theorem.

\begin{thm}[Inverse Function Theorem]
Let $F:\mathcal{B}\too\mathcal{B}'$ be a $C^1$ map between Banach
spaces and suppose that the linearization $\DD F(0)$ is an
isomorphism.  Moreover, assume that for some constants $C_I$,
$C_N$, and $r_1$, we have
\begin{quote}
\begin{itemize}
\item[1.] $\|\DD F(0)x\|_{\mathcal{B}'}\geq
  {1\over C_I}\|x\|_{\mathcal{B}}$ for all
  $x\in\mathcal{B}$, and
\item[2.] $\|\DD F(0)x-\DD F(y)x\|_{\mathcal{B}'}\leq
 C_N\|x\|_{\mathcal{B}}\!\cdot\!\|y\|_{\mathcal{B}}$
  for all $x,y\in\mathcal{B}$ with $\|y\|_{\mathcal{B}}<r_1$.
\end{itemize}
\end{quote}
 Then there exist neighborhoods $U$ of $0$ and $V$ of $F(0)$ such that
 $F:U\too V$ is a $C^1$-diffeomorphism. Moreover,
 if $r\leq \mathrm{min}(r_1,(2 C_I C_N)^{-1})$, then
 $B_{r/2C_I}(F(0))\subset V$ and $B_{r/2C_I}(F(0))\subset
 F(B_r(0))$.
\end{thm}
In particular, when the hypotheses of the theorem are satisfied
and additionally, $\|F(0)\|_{\mathcal{B}'}< {r\over 2C_I}$, we can
solve the equation $F(y)=0$ for some $\|y\|_{\mathcal{B}}<r$.

In order to invoke the Inverse Function Theorem in our situation,
we need to choose our Banach spaces carefully.  We define a smooth
weight function $\rho$ on $\M$ with the key property that the ball
of radius $\rho(x)$ in $\M$ centered at $x$ has uniformly bounded
geometry.  That is, in geodesic normal coordinates at $x$, we have
$|g_{ij}-\delta_{ij}|^*_{1,\beta,B_{\rho(x)}(x)}\leq 1$ where the
norm here is the local scale-invariant Schauder norm on
$B_{\rho(x)}(x)$. We also require that $\tau=1$ on the ball
$B_{\rho(x)}(x, X)$.  We can construct such a $\rho$ with the
following additional properties. See
\cite{butscher,butscherthesis,ylee}.
\begin{quote}
\begin{itemize}
\item[$\bullet$] For some $r_0$ and $R$ independent of $\alpha$,
$$\rho(x)=\left\{
\begin{array}{ll}
\epsilon R &\text{ for }x\in N'=\M\cap B_{\epsilon r_0}\\\
\text{interpolation} &\text{ for }x\in \M\cap(B-B_{\epsilon r_0})\\
R &\text{ for }x\in \M-B
\end{array}\right.$$

\item[$\bullet$] $\rho(x)\leq C|x|$ for $x\in
\M\cap(B-B_{\delta/2})$. \item[$\bullet$] $|\nabla\rho|_0\leq C$.
\item[$\bullet$]$\|\rho^{-1}\|_{L^2(\M)}\leq C$.
\end{itemize}
\end{quote}

\begin{defin}
For any $0<\beta<1$, the \textbf{$\rho$-weighted
$(k,\beta)$-Schauder norm} on $C^{k,\beta}(\M)$ is given by
$$|u|_{C^{k,\beta}_\rho(\M)}=|u|_{0,\M}+|\rho\nabla
u|_{0,\M}+\cdots+|\rho^k\nabla^k u|_{0,\M}+[\rho^{k+\beta}\nabla^k
u]_{\beta,\M}.$$ Let $S$ be the first eigenfunction of the
Laplacian on $\M$, normalized so that
$\|S\|_{L^2(\M)}=1$.\footnote{In contrast to our use of constants,
many geometric objects such as functions and operators will depend
on $\alpha$, but we will suppress this dependence in the notation
for the purpose of readability.  The loss of clarity should be
minimal since these objects are all defined on $\M$.} Then we
define the Banach spaces $\mathcal{B}_{1,\alpha}$,
$\mathcal{B}_\alpha$, and $\mathcal{B}'_\alpha$ as vector spaces
\begin{eqnarray*}
\mathcal{B}_{1,\alpha}&=&\left\{u\in C^{2,\beta}(\M)\left|
\int_{\M} u =\int_{\M} u S = 0\right.\right\}\\
\mathcal{B}_{\alpha}&=&\mathcal{B}_{1,\alpha}\times\rr\\
\mathcal{B}'_\alpha&=&\left\{u\in C^{0,\beta}(\M)\left| \int_{\M}
u=0\right.\right\}
\end{eqnarray*}
with the norms
\begin{eqnarray*}
\|u\|_{\mathcal{B}_{1,\alpha}}&=&|u|_{C^{2,\beta}_\rho(\M)}\\
\|(u,a)\|_{\mathcal{B}_{\alpha}}&=&|u|_{C^{2,\beta}_\rho(\M)}+|a|\\
\|f\|_{\mathcal{B}'_\alpha}&=&|\rho^2 f|_{C^{0,\beta}_\rho(\M)}.
\end{eqnarray*}
(The integrations above are taken with respect to the
$t$-independent K\"ahler metric on $\M$.)
\end{defin}
Since $\Psi_t$ and $\Phi_h$ are isotopies,
$\langle[(\Psi_t\circ\Phi_h)(\M)],[\im\Omega_t]\rangle=
\langle[\M],[\im\Omega_t]\rangle= \langle
[M_1]+[M_2],[\im\Omega_t]\rangle=0$, and therefore
$F_\alpha(\mathcal{B}_{\alpha})\subset\mathcal{B}'_{\alpha}$. From
now on we think of the \textbf{deformation operator} $F_\alpha$ as
an operator from $\mathcal{B}_\alpha$ to $\mathcal{B}'_\alpha$.

The choice of $\beta$ is not particularly important; it is simply
a small constant independent of $\alpha$. The purpose of the
weighted norm is to achieve estimates that scale nicely with
respect to $\alpha$. The reason why we take the orthogonal
complement of the functions $1$ and $S$ is that $1$ lies in the
kernel of the linearization, and $S$ lies in the approximate
kernel of the linearization.

Let us summarize what we need to prove in order to invoke the
Inverse Function Theorem argument:
\begin{quote}
\begin{itemize}
\item[$\bullet$] We need an injectivity estimate on $\DD
F_\alpha(0,0)$; we must establish the existence of a constant
$C_I(\alpha)$ as in condition 1 of the Inverse Function Theorem
and find its dependence on $\alpha$. \item[$\bullet$] We need to
show that $\DD F_\alpha(0,0)$ is surjective. \item[$\bullet$] We
need a nonlinear estimate; we must establish the existence of a
constant $C_N(\alpha)$ as in condition 2 of the Inverse Function
Theorem and find its dependence on $\alpha$. \item[$\bullet$] We
need to bound $F_\alpha(0,0)$ in terms of $\alpha$.
\end{itemize}
\end{quote}
We first compute $\DD F_\alpha(0,0)$.
\begin{prop}
\begin{equation}\label{DF}
\DD F_\alpha(0,0)(u,a)=\Delta u + Pu+a\psi
\end{equation}
where $P:\mathcal{B}_{1,\alpha}\too\mathcal{B}'_\alpha$ is given
by
$$
Pu=(\cos\theta-1)\Delta u-(\sin\theta)\langle H,J\nabla u\rangle
$$
and $\psi\in\mathcal{B}'_\alpha$ is given by
$$\psi=\langle L_V(\im\Omega)+\im\dot{\Omega},\vol_{\M}\rangle$$
where $V$ is the vector field generating the flow $\Psi_t$ at time
$t=0$, and \mbox{$\dot{\Omega}={\dd\over \dd t}\Omega_t|_{t=0}$}.
\end{prop}
The calculation of $\psi$ is self-evident. The rest of the
calculation is straightforward and can be found in \cite{butscher}
and \cite{ylee}.  The reason we write $\Delta$ and $P$ separately
in equation (\ref{DF}) is that the $P$ term turns out to be
negligible, and therefore it suffices to understand $\Delta$ and
$\psi$. The unimportance of $P$ is expressed in the following
lemma.
\begin{lem}\label{Pu}
For sufficiently small $\alpha$, for any
$u\in\mathcal{B}_{1,\alpha}$,
$$\|Pu\|_{\mathcal{B}'_\alpha}\leq C \alpha^{1-\beta}
\|u\|_{\mathcal{B}_{1,\alpha}}.$$
\end{lem}
\begin{proof}
The proof essentially follows directly from the bounds given in
\mbox{Lemma \ref{theta}}.
\begin{eqnarray*}
|\rho^2(1-\cos\theta)\Delta u|_0 &\leq& |1-\cos\theta|_0\!\cdot\!
|\rho^2\Delta u|_0\\
&\leq&C\alpha^2 |u|_{C_\rho^{2,\beta}}.
\end{eqnarray*}
\begin{eqnarray*}
[\rho^{2+\beta}(1-\cos\theta)\Delta u]_\beta &\leq&
[1-\cos\theta]_\beta\!\cdot\! |\rho^{2+\beta}\Delta u|_0 +
|1-\cos\theta|_0\!\cdot\! [\rho^{2+\beta}\Delta u]_\beta\\
&\leq&(C\alpha^{2-\beta}R^\beta+C\alpha^2)|u|_{C_\rho^{2,\beta}}.
\end{eqnarray*}
\begin{eqnarray*}
|\rho^2(\sin\theta)\langle H,J\nabla u\rangle|_0 &\leq&
|\rho(\sin\theta)H|_0\!\cdot\!|\rho J\nabla u|_0\\
&\leq&RC\alpha |u|_{C_\rho^{2,\beta}}.
\end{eqnarray*}
For the final inequality, we use the fact that $[v]_\beta\leq
|v|_0+|\nabla v|_0$.
\begin{eqnarray*}
[\rho^{2+\beta}(\sin\theta)\langle H,J\nabla u\rangle]_\beta
&\leq& [\sin\theta]_\beta\!\cdot\!|H|_0\!\cdot\!|\rho^{2+\beta}
J\nabla u|_0+
|\sin\theta|_0\!\cdot\![H]_\beta\!\cdot\!|\rho^{2+\beta} J\nabla
u|_0 +\\
& & |\sin\theta|_0\!\cdot\!|H|_0\!\cdot\![\rho^{2+\beta}
J\nabla u]_\beta\\
&\leq& (C\alpha^{1-\beta})CR^{1+\beta}|\rho\nabla u|_0 +
(C\alpha)(C\alpha^{-\beta})R^{1+\beta}|\rho\nabla u|_0 +\\
& &(C\alpha)C(|\rho^{2+\beta} \nabla u|_0+|(\nabla\rho^{2+\beta})
(\nabla u)|_0+|\rho^{2+\beta}\nabla^2 u|_0)\\
&\leq&C\alpha^{1-\beta}|u|_{C_\rho^{2,\beta}} +C\alpha
(R^{1+\beta}+CR^\beta  + R^\beta)|u|_{C_\rho^{2,\beta}}
\end{eqnarray*}
where the last line uses the bound on $|\nabla\rho|_0$.  Now
combine the previous four inequalities to deduce the desired
result.
\end{proof}

\section{Analysis of the Laplacian on $\M$}

The first step in establishing an injectivity estimate for the
linearized deformation operator is finding a lower bound for the
second eigenvalue of the Laplacian.  The second step is to combine
this lower bound with an elliptic estimate to obtain an
injectivity estimate for
$\Delta:\mathcal{B}_{1,\alpha}\too\mathcal{B}'_\alpha$.

It is a fact that on any Riemannian manifold $M$, for any $f\in
L^2(M)$ with one derivative in $L^2(M)$ such that $\int_M f=0$, we
have $\int_M |\nabla f|^2\geq \lambda_1(M)\int_M f^2$, where
$\lambda_1(M)$ is the first eigenvalue of the Laplacian.  From
this it follows easily that if we drop the condition $\int_M f=0$,
then we have
\begin{equation}\label{evalest}
\lambda_1(M)\leq {\int_M |\nabla f|^2\over\int_M f^2-{1\over
\vol(M)}(\int_M f)^2}.
\end{equation}

We define a smooth cutoff function $\varphi$ on $\M$ with the
following properties: $\varphi=0$ in $B_\delta$, $\varphi=1$
outside $B_{2\delta}$, and $|\nabla\varphi|\leq {C \over \delta}$
for some $C$.  Recall that $\M-B_\delta=(M_1\cup M_2)-B_\delta$,
and therefore we may think of $\varphi$ as a function on either
$\M$ or on $M_1\amalg M_2$.  Observe that because we have
uniformly bounded mean curvature, the Monotonicity Formula
provides the following bounds which we will use repeatedly:
\begin{eqnarray*}
\vol(\M\cap B_{2\delta})&=&O(\delta^n)\\
\vol(M_i\cap B_{2\delta})&=&O(\delta^n)\\
\vol(M_i-B_{2\delta})&\leq&C\\
\vol(M_i-B_{2\delta})&\geq&{1\over C}
\end{eqnarray*}

\begin{lem}\label{firsteval}
For small enough $\alpha$,
  $$\lambda_1(\M)\leq C\delta^{n-2}.$$
\end{lem}
\begin{proof}
Let
$$f=\varphi\left[{\mathcal{X}_{M_1-B_\delta}\over
\vol(M_1-B_{2\delta})}-{\mathcal{X}_{M_2-B_\delta}\over
\vol(M_2-B_{2\delta})}\right]$$ where $\mathcal{X}_A$ denotes the
characteristic function of $A$.  Clearly, $\int_{\M}
f=O(\delta^n)$.  Using inequality (\ref{evalest}), we see that
\begin{eqnarray*}
\lambda_1(\M) &\leq& {\int_{\M} |\nabla f|^2\over
 \int_{\M} f^2-{1\over \vol(\M)}(\int_{\M}f)^2 }\\
&\leq& {C\int_{\M\cap B_{2\delta}}\delta^{-2}\over {1\over
\vol(M_1-B_{2\delta})}+{1\over \vol(M_2-B_{2\delta})}
-O(\delta^{2n})}\text{ by the properties of }\varphi\\
&\leq&C\delta^{n-2}\text{ since the denominator is bounded below.}
\end{eqnarray*}
\end{proof}
We would like to have some idea of what $S$ looks like.  By the
previous Lemma together with Lemma 5 of Yng-Ing Lee's paper
\cite{ylee}, we know that $|S|_0$ is bounded independently of
$\alpha$.\footnote{This Lemma depends on the Michael-Simon
Inequality and uses the fact that the mean curvature of $\M$ is
bounded independently of $\alpha$.}  This fact allows us to use
our knowledge of the kernel of the Laplacian on $M_1\amalg M_2$ to
construct a function that approximates $S$ in the $L^2$ sense.
\begin{lem}\label{sbar}
Define
$$\bar{S}=a_1\mathcal{X}_{M_1}+a_2\mathcal{X}_{M_2}$$ where
$$a_1={1\over \vol(M_1)}
\sqrt{{\vol(M_1)\vol(M_2)\over \vol(M_1)+\vol(M_2)}}\text{ and }
a_2={-1\over \vol(M_2)} \sqrt{{\vol(M_1)\vol(M_2)\over
\vol(M_1)+\vol(M_2)}}.$$
 Then for small enough $\alpha$,
$$\|S-\varphi\bar{S}\|_{L^2(\M)}\leq C\delta^{(n-2)/2}$$
and
$$\|\bar{S}-\varphi\bar{S}\|_{L^2(M_1\amalg M_2)}\leq C\delta^{n/2}.$$
\end{lem}
\begin{proof}
First, the bound on $S$ implies that
\begin{equation}\label{sbar1}
\|S-\varphi S\|^2_{L^2(\M)}\leq \|S\|^2_{L^2(\M\cap
B_{2\delta})}=O(\delta^{n}).
\end{equation}
 Since $\varphi S$ is defined on
$M_1\amalg M_2$, we have
\begin{eqnarray*}
\int_{M_1\amalg M_2} |\nabla(\varphi S)|^2  &=&\int_{\M}
|\nabla(\varphi S)|^2\\
&\leq&2\int_{\M\cap B_{2\delta}}|\nabla\varphi|^2 S^2
+2\int_{\M}\varphi^2 |\nabla S|^2  \\
&\leq&O(\delta^{n-2})+2\int_{\M}|\nabla S|^2 \text{ arguing as in
Lemma \ref{firsteval}}\\
&=&O(\delta^{n-2})\text{ by Lemma \ref{firsteval} and the
normalization of }S.
\end{eqnarray*}
Note that $M_1\amalg M_2$ has a two-dimensional kernel spanned by
$\mathcal{X}_{M_1}$ and $\mathcal{X}_{M_2}$, and its first
non-zero eigenvalue is obviously a constant independent of
$\alpha$. Therefore the estimate above shows that if
$a'_1\mathcal{X}_{M_1}+a'_2\mathcal{X}_{M_2}$ is the orthogonal
projection of $\varphi S$ onto the kernel, then
\begin{eqnarray}
\|\varphi S -
(a'_1\mathcal{X}_{M_1}+a'_2\mathcal{X}_{M_2})\|^2_{L^2(M_1\amalg
M_2)}&\leq& C\|\nabla[\varphi S -
(a'_1\mathcal{X}_{M_1}+a'_2\mathcal{X}_{M_2})]\|^2_{L^2(M_1\amalg
M_2)}\nonumber\\
&=&C\|\nabla(\varphi S)\|^2_{L^2(M_1\amalg
M_2)}\nonumber\\
&=&O(\delta^{n-2})\text{ by the previous
calculation.}\label{sbar2}
\end{eqnarray}
The bounds (\ref{sbar1}) and (\ref{sbar2}) show that
$$(a'_1)^2\vol(M_1)+(a'_2)^2\vol(M_2)=
\|a'_1\mathcal{X}_{M_1}+a'_2\mathcal{X}_{M_2}\|^2_{L^2(M_1\amalg
M_2)}=1+ O(\delta^{n-2})$$ and
$$a'_1 \vol(M_1)+a'_2\vol(M_2) =
\int_{M_1\amalg M_2}a'_1\mathcal{X}_{M_1}+a'_2\mathcal{X}_{M_2} =
O(\delta^{(n-2)/2}).$$ Solving these equations, we find that
$$a'_1=a_1+O(\delta^{(n-2)/2})\text{ and }
a'_2=a_2+O(\delta^{(n-2)/2})$$ where $a_1$ and $a_2$ were defined
above.\footnote{Of course, the solution is only determined up to a
sign, but there was a sign ambiguity in our original definition of
$S$, so we can simply define $S$ to have the sign consistent with
these equations.}
 It now follows that
\begin{equation}\label{sbar3}
\|(a'_1\mathcal{X}_{M_1}+a'_2\mathcal{X}_{M_2})-
\bar{S}\|_{L^2(M_1\amalg M_2)}^2=O(\delta^{n-2}).
\end{equation}
Finally, similar to the bound (\ref{sbar1}), we see that
\begin{equation}\label{sbar4}
\|\bar{S}-\varphi\bar{S}\|^2_{L^2(M_1\amalg M_2)}=O(\delta^{n}).
\end{equation}
 Putting together the bounds
(\ref{sbar1}), (\ref{sbar2}), (\ref{sbar3}), and(\ref{sbar4}), we
obtain the desired result.
\end{proof}
We can now use our knowledge of $S$ to help us show that the
second eigenvalue of the Laplacian on $\M$ is bounded below.
\begin{prop}[Second Eigenvalue Estimate]
For small enough $\alpha$, the second eigenvalue of the Laplacian
of $\M$, $\lambda_2(\M)$, is bounded below.  In particular, for
each $u\in\mathcal{B}_{1,\alpha}$,
$$\|\Delta u\|_{L^2(\M)}\geq {1\over C}\|u\|_{L^2(\M)}$$
\end{prop}
\begin{proof}
Let $f$ be an eigenfunction for the second eigenvalue of the
Laplacian, normalized so that $\|f\|_{L^2(\M)}=1$.  Using the
min-max characterization of $\lambda_2(\M)$ one can show that
$\lambda_2(\M)$ is bounded above independently of $\alpha$. This
fact allows us to apply Lemma 5 of Yng-Ing Lee's paper \cite{ylee}
to show that $|f|_0$ is bounded independently of $\alpha$. We
compute

\begin{eqnarray*}
\lambda_2(\M)&=& \int_{\M}|\nabla f|^2\\
&=& \int_{\M}|\nabla (\varphi f)+\nabla
((1-\varphi)f)|^2\\
&\geq&\int_{M_1}|\nabla (\varphi f)|^2+\int_{M_2}|\nabla (\varphi
f)|^2- 2\int_{\M\cap B_{2\delta}}|\nabla(1-\varphi)|^2f^2-
2\int_{\M}(1-\varphi)^2|\nabla f|^2)\\
&\geq&\int_{M_1}|\nabla (\varphi f)|^2+\int_{M_2}|\nabla (\varphi
f)|^2- O(\delta^{n-2}) -2\int_{\M}|\nabla f|^2.
\end{eqnarray*}
Therefore,
\begin{eqnarray*}
3\lambda_2(\M) &\geq& \int_{M_1}|\nabla (\varphi
f)|^2+\int_{M_2}|\nabla (\varphi f)|^2- O(\delta^{n-2})\\
&\geq&\lambda_1(M_1)\left[\int_{M_1}\varphi^2 f^2 -{1\over
\vol(M_1)}\left(\int_{M_1}\varphi f\right)^2\right] +\\
& &\lambda_2(M_2)\left[\int_{M_2}\varphi^2 f^2 -{1\over
\vol(M_2)}\left(\int_{M_2}\varphi f\right)^2\right] -
O(\delta^{n-2})\text{ by inequality (\ref{evalest})}\\
&\geq&C\int_{\M} f^2-O(\delta^n) -{\lambda_1(M_1)\over
\vol(M_1)}\left(\int_{M_1}\varphi f\right)^2\\
& & - {\lambda_2(M_2)\over \vol(M_2)}\left(\int_{M_2}\varphi
f\right)^2 -
O(\delta^{n-2})\\
&\geq&C-C\left[\left(\int_{M_1}\varphi
f\right)^2+\left(\int_{M_2}\varphi
f\right)^2\right]-O(\delta^{n-2}).
\end{eqnarray*}
Now we must bound the terms in brackets in the previous line.
\begin{eqnarray*}
\left|\int_{M_1}\varphi f\right| &=& \left|\int_{M_1\amalg
M_2}\left({\bar{S}-a_2\over a_1-a_2}\right)\varphi f\right|\text{
by
the definition of }\bar{S}\text{ from Lemma \ref{sbar}}\\
&=&{1\over
a_1-a_2}\left|\int_{\M}(\varphi\bar{S})f-a_2\int_{\M}\varphi
f\right|\\
&=&{1\over a_1-a_2}\left|\int_{\M}Sf+\int_{\M}(\varphi\bar{S}-S)f
-a_2\int_{\M}f+O(\delta^n)\right|\\
&\leq&C\|\varphi\bar{S}-S\|_{L^2(\M)}+ O(\delta^n) \text{ since
}\int_{\M}Sf=\int_{\M}f=0\\
&=&O(\delta^{(n-2)/2})+O(\delta^n)\text{ by Lemma \ref{sbar}.}
\end{eqnarray*}
The estimate for $M_2$ is similar.
\end{proof}

Recall that the weight function $\rho$ was chosen so that $\M$ has
uniformly bounded geometry in $\rho(x)$ neighborhoods of $x$ in
$\M$.  This allows us to use the local scale-invariant elliptic
Schauder estimate to deduce a global elliptic Schauder estimate
independent of $\alpha$.  We omit the proof which is standard and
straightforward.
\begin{prop}[Global Elliptic Schauder Estimate]
For sufficiently small $\alpha$, for any $u\in C^{2,\beta}(\M)$,
$$|u|_{C_\rho^{2,\beta}}\leq C(|\rho^2\Delta
u|_{C_\rho^{0,\beta}}+|u|_0).$$
\end{prop}
The following lemma, proved in \cite{ylee}, translates the second
eigenvalue estimate from the $L^2$ setting to the Schauder
setting.\footnote{The proof of this lemma follows easily from a De
Giorgi-Nash estimate, which in turn depends on the Michael-Simon
inequality and bounded mean curvature.}  Choose any $\nu>0$
independent of $\alpha$.
\begin{lem}
For sufficiently small $\alpha$,
$$|u|_0\leq \epsilon^{-\nu}|\rho^2\Delta u|_{C_\rho^{0,\beta}}$$
for all $u\in\mathcal{B}_{1,\alpha}$.
\end{lem}
Combining this Lemma with Global Elliptic Estimate immediately
leads us to an injectivity estimate for
$\Delta:\mathcal{B}_{1,\alpha}\too\mathcal{B}'_\alpha$.
\begin{prop}[Laplacian Injectivity Estimate]
For sufficiently small $\alpha$, for all
$u\in\mathcal{B}_{1,\alpha}$,
$$|u|_{C_\rho^{2,\beta}}\leq C\epsilon^{-\nu}|\rho^2\Delta
u|_{C_\rho^{0,\beta}}.$$
\end{prop}

\section{Proofs of the Key Lemma}\label{proofkey}

It is well-known that $\Delta$ is an isomorphism from
$\mathcal{B}_{1,\alpha}\oplus\langle S\rangle$ to
$\mathcal{B}'_\alpha$, but because of the small first eigenvalue,
we had to remove $S$ in order to obtain a good injectivity
estimate.  Of course, removing $S$ costs us surjectivity.  We
added the extra degree of freedom in order to restore
surjectivity.  The essential requirement of the extra degree of
freedom is that its linearization $\psi$ must have a significant
$S$ component.  This is the content of our Key Lemma.  We now
construct a deformation $(J_t,\omega_t,\Omega_t)$ so that $\psi$
has the desired property.

\begin{lem}\label{deform}
Under the assumptions of the Main Theorem or the Torus Version,
there exists a deformation of Calabi-Yau structure
$(J_t,\omega_t,\Omega_t)$ such that $\omega_t$ is cohomologous to
$\omega$, $\langle [M_1]+[M_2], [\im\Omega_t]\rangle=0$, and
$$\left\langle {[M_1]\over\vol(M_1)}-{[M_2]\over \vol(M_2)},
[\im\dot{\Omega}]\right\rangle\neq 0$$ for sufficiently small
$\alpha$.
\end{lem}
\begin{proof}[Proof of Main Case]
By a result of Tian \cite{tian} and Todorov \cite{todorov}, the
first-order deformation space of complex structures on $(X,J)$ is
exactly the space of harmonic $(n-1,1)$ forms,
$\mathcal{H}^{n-1,1}(X)$, and all of these first-order
deformations extend to actual deformations.  By the hypotheses of
the Main Theorem, we can find
$\lambda\in\re[\mathcal{H}^{n-1,1}(X)\oplus\mathcal{H}^{1,n-1}(X)]$
such that $\left\langle {[M_1]\over\vol(M_1)}-{[M_2]\over
\vol(M_2)}, [\lambda]\right\rangle\neq0$.  Since
$\mathcal{H}^{n-1,1}(X)$ and $\mathcal{H}^{1,n-1}(X)$  are complex
conjugate to each other, we can certainly find
$\chi\in\mathcal{H}^{n-1,1}(X)$ such that $\im\chi=\lambda$.
Choose a complex structure deformation $J_t$ whose first-order
deformation is $\chi$.  Recall that the choice of $J_t$ determines
the holomorphic $(n,0)$-form $\Omega_t$ up to a constant.  It is a
fact that the effect of the first-order deformation $\chi$ on
$\dot{\Omega}$ is expressed by the formula
$$\dot{\Omega}=c\Omega+\chi.$$
See Candelas and de la Ossa \cite{canoss} for details. Since
$\im\Omega$ vanishes on $M_1$ and $M_2$, and $\re\Omega$
calibrates $M_1$ and $M_2$, it follows that $\left\langle
{[M_1]\over\vol(M_1)}-{[M_2]\over \vol(M_2)},
[\Omega]\right\rangle=0$, and therefore $\left\langle
{[M_1]\over\vol(M_1)}-{[M_2]\over \vol(M_2)},
[\im\dot{\Omega}]\right\rangle\neq 0$.

We know that $[\omega]$ lies in the K\"ahler cone of $(X,J)$.
Since $(X,J)$ admits a Calabi-Yau structure with holonomy equal to
$\su(n)$, so does $(X,J_t)$.  In this case, it is a fact that
$H^{2,0}(X,J_t)=H^{0,2}(X,J_t)=0$. See Joyce \cite{joycebook} for
details.  Therefore the K\"ahler cone of $(X,J_t)$ is open in
$H^2(X)$ for each $t$.  So for small $t$, $[\omega]$ lies in the
K\"ahler cone of $(X,J_t)$.  Therefore we can find a unique
$\omega_t\in[\omega]$ such that $(X,J_t,\omega_t)$ admits a CY
structure.  This determines $\Omega_t$ up to a phase, and we
choose the phase so that $\langle [M_1]+[M_2],
[\im\Omega_t]\rangle=0$.
\end{proof}

Before we discuss the Torus Case, let us study Calabi-Yau tori and
the moduli space of CY tori.
\begin{prop}\label{lattice}
We can characterize all Calabi-Yau tori and all flat special
Lagrangian submanifolds of them as follows.
\begin{quote}
\begin{itemize}
\item[1.] Given a rank $2n$ lattice $\Gamma\subset\cc^n$, we
define a Calabi-Yau structure on $\cc^n/\Gamma$ using the standard
Calabi-Yau structure on $\cc^n$. Every Calabi-Yau torus $T$ is
given by this construction. \item[2.] Given a special Lagrangian
plane $\eta$ invariant with respect to some rank $n$ sublattice of
$\Gamma$, we obtain a special Lagrangian torus $\eta/\Gamma$ in
$\cc^n/\Gamma$. Every flat special Lagrangian submanifold of $T$
is a union of special Lagrangian tori given by this construction.
\end{itemize}
\end{quote}
\end{prop}
\begin{proof}
Let $(T,J,\omega,\Omega)$ be a CY torus.  Clearly, we have an
induced CY structure on the universal cover
$(\rr^{2n},\tilde{J},\tilde{\omega},\tilde{\Omega})$.  Since the
K\"ahler metric on $T$ has zero Ricci curvature, we can apply
Cheeger and Gromoll's theorem on manifolds with non-negative Ricci
curvature \cite{chegro} to see that the induced metric on the
universal cover is flat.  Choose a point $q\in\rr^{2n}$.  We can
find a basis $e_1,\ldots,e_{2n}$ of $T_q \rr^{2n}$ such that if
$e^1,\ldots,e^{2n}$ is the dual basis, then $\tilde{J}e_j=e_{j+n}$
for $j\leq n$, $\tilde{J}e_j=-e_{j-n}$ for $j>n$,
$\tilde{\omega}=\sum_{j=1}^n e^j\wedge e^{j+n}$, and
$\tilde{\Omega}=(e_1+ie_{n+1})\wedge\ldots\wedge(e_n+ie_{2n})$ at
the point $q$.  Since the metric is flat, we can extend
$e_1,\ldots,e_{2n}$ to a frame field over all of $\rr^{2n}$ in a
unique way via parallel translation.  Note that $[e_j,e_k]=0$
everywhere, so we can find a global coordinate system
$x^1,\ldots,x^n,y^1,\ldots,y^n$ such that ${\partial\over\partial
x^j}=e_j$ and ${\partial\over\partial y^j}=e_{j+n}$.  Since
$\nabla J=\nabla \omega=\nabla \Omega=0$, it follows that
$\tilde{J}({\partial\over\partial x^j})={\partial\over\partial
y^j}$, $\tilde{J}({\partial\over\partial
y^j})=-{\partial\over\partial x^j}$, $\tilde{\omega}=\sum_{j=1}^n
\dd x^j\wedge \dd y^j$, and $\tilde{\Omega}=\dd
z^1\wedge\ldots\wedge \dd z^n$ everywhere.  That is, the induced
CY structure on the universal cover is the standard CY structure
on $\cc^n$.  Moreover, since each Deck transformation must be an
orientation-preserving isometry, the Deck transformations are
translations.

Let $M$ be a flat SLag in $T$.  Since $M$ is flat and minimal in
$T$, which is also flat, it is a simple consequence of the Gauss
equation that $M$ is totally geodesic in $T$.  In fact,
Ricci-flatness would have been sufficient.  See \cite{dajczer}.
Therefore, locally, the lift of $M$ up to $\cc^n$ is a piece of a
SLag plane. It follows that $M$ is the quotient of a union of SLag
planes which are invariant under some rank $n$ sublattice of
$\Gamma$.
\end{proof}
We know that the CY structure of a torus is determined by a
lattice.  The space of all rank $2n$ lattices in $\cc^n$ is
$\gl(2n,\rr)/\slin(2n,\zz)$ where $\slin(2n,\zz)$ acts on the
right. Since the group of CY structure preserving automorphisms of
$\cc^n$ is $\su(n)\rtimes$(translations), it follows that the
global moduli space of Calabi-Yau tori is precisely
$\gl(2n,\rr)/\slin(2n,\zz)$ modulo the action of $\su(n)$ on the
left. Therefore, locally, the deformation space is simply a
neighborhood of the identity in $\gl(2n,\rr)/\su(n)$ where
$\su(n)$ acts on the left. Ignoring the (discrete) redundancies
arising from the lattice automorphisms, $\slin(2n,\zz)$, the space
of possible complex structures corresponds to
$\gl(2n,\rr)/\gl(n,\cc)$, and after this choice is made, the space
of compatible symplectic structures corresponds to
$\gl(n,\cc)/\un(n)$.  On the other hand, we can choose the
symplectic structure from the space $\gl(2n,\rr)/\symp(2n,\rr)$
and then choose a compatible complex structure from the space
$\symp(2n,\rr)/\un(n)$.  Finally, of course, we choose a phase for
$\Omega$ from $\un(n)/\su(n)$. Finally, observe that changing the
lattice by $A\in\gl(2n,\rr)$ has the same effect as keeping the
lattice and the canonical local coordinate systems on $T$ fixed,
but changing the CY structure to $(A^{-1}JA,A^*\omega,A^*\Omega)$
with respect to those coordinates.  This is the point of view we
adopt in the following proof, and we will no longer mention
lattices.

\begin{proof}[Proof of Torus Case of Lemma \ref{deform}]
By Proposition \ref{lattice}, any connected, flat SLag $M_1$ in
$(T, J,\omega,\Omega)$ is actually a SLag torus with constant
tangent plane with respect to the canonical local coordinates, so
we can perform an $\su(n)$ change of coordinates taking $M_1$ to a
SLag torus with tangent plane ${\partial\over\partial
x^1}\wedge\ldots\wedge{\partial\over\partial x^n}$ at each point
of $M_1$.  Since the $\su(n)$ change of coordinates preserves the
CY structure, we may assume without loss of generality that $T_q
M_1= {\partial\over\partial
x^1}\wedge\ldots\wedge{\partial\over\partial x^n}$ at each $q\in
M_1$.

Consider a deformation $A_t^{-1}JA_t$ of the complex structure.
Then the holomorphic $(n,0)$ form is $A_t^*\Omega$ up to a
constant.  A simple calculation then shows that
$\im\dot{\Omega}=c\Omega+\chi$, where
$$\chi=\dsum_{j,k=1}^n B_{jk}\dd
z^1\wedge\ldots\wedge \overbrace{\dd\bar{z}^j}^{k\text{-th
spot}}\wedge\ldots\wedge \dd z^n$$ and
$B_{jk}=\dot{A}_{k\bar{j}}$. Another simple calculation shows that
if $A_t\in\symp(2n,\rr)$, then the corresponding $B$ must be
(complex) symmetric, and conversely, for any (complex) symmetric
$B$, we can find $A_t\in\symp(2n,\rr)$ such that
$B_{jk}=\dot{A}_{k\bar{j}}$.

In particular, we can find $A_t\in\symp(2n,\rr)$ and a phase
$\phi_t\in\rr$ such that
$(A_t^{-1}JA_t,\omega,\Omega_t=e^{i\phi_t}A_t^*\Omega)$ is a CY
structure with $\langle [M_1]+[M_2],[\im\Omega_t]\rangle=0$, and
$$\chi=i\dsum_{j=1}^n \dd
z^1\wedge\ldots\wedge \overbrace{\dd\bar{z}^j}^{j\text{-th
spot}}\wedge\ldots\wedge \dd z^n.$$ Note that
$$\im\chi=
\dsum_{j=1}^n \re(\dd z^1\wedge\ldots\wedge
\overbrace{\dd\bar{z}^j}^{j\text{-th spot}}\wedge\ldots\wedge \dd
z^n).$$  It remains to show that $\left\langle
{[M_1]\over\vol(M_1)}-{[M_2]\over \vol(M_2)},
[\im\chi]\right\rangle\neq0$. First, observe that $\re(\dd
z^1\wedge\ldots\wedge \dd\bar{z}^j\wedge\ldots\wedge \dd z^n)$ is
a calibration for each $j$, and that ${\partial\over\partial
x^1}\wedge\ldots\wedge{\partial\over\partial x^n}$ is calibrated
by each of these calibrations.  Therefore
${1\over\vol(M_1)}\int_{M_1}\im\chi=n$ and
${1\over\vol(M_2)}\int_{M_2}\im\chi\leq n$ with equality iff $M_2$
is also calibrated by each of the calibrations $\re(\dd
z^1\wedge\ldots\wedge \dd\bar{z}^j\wedge\ldots\wedge \dd z^n)$. We
momentarily consider the ``torus" planes $\xi$ of the form
$$\left[[(\cos\theta_1){\partial\over\partial x^1}+
(\sin\theta_1){\partial\over\partial y^1}\right]\wedge\ldots\wedge
\left[[(\cos\theta_n){\partial\over\partial x^n}+
(\sin\theta_n){\partial\over\partial y^n}\right].$$ It is easy to
verify that if $\xi$ is calibrated by $\re\Omega=\re(\dd
z^1\wedge\ldots\wedge \dd z^n)$ and $\re(\dd z^1\wedge\ldots\wedge
\dd\bar{z}^j\wedge\ldots\wedge \dd z^n)$ for each $j$, then
$\xi={\partial\over\partial
x^1}\wedge\ldots\wedge{\partial\over\partial x^n}$.  Then by
Morgan's Torus Lemma, it follows that ${\partial\over\partial
x^1}\wedge\ldots\wedge{\partial\over\partial x^n}$ is the only
plane simultaneously calibrated by $\re(\dd z^1\wedge\ldots\wedge
\dd z^n)$ and $\re(\dd z^1\wedge\ldots\wedge
\dd\bar{z}^j\wedge\ldots\wedge \dd z^n)$ for each $j$.  See
\cite{harvey} for more on Morgan's Torus Lemma. Therefore if $M_2$
is calibrated by $\re\Omega=\re(\dd z^1\wedge\ldots\wedge \dd
z^n)$ and $\re(\dd z^1\wedge\ldots\wedge
\dd\bar{z}^j\wedge\ldots\wedge \dd z^n)$ for each $j$, it follows
that $T_p M_2={\partial\over\partial
x^1}\wedge\ldots\wedge{\partial\over\partial x^n}$, violating the
transversality assumption.
\end{proof}

\begin{lem}[Key Lemma]
For small $\alpha$,
\begin{equation}
\left|\int_{\M}\psi S\right|\geq {1\over C}.
\end{equation}
\end{lem}
\begin{proof}
Recall that $\psi=\langle L_V(\im\Omega),\vol_{\M}\rangle
+\langle\im\dot{\Omega},\vol_{\M}\rangle$.  We'll show that the
first term integrated against $S$ is small while the second term
integrated against $S$ is bounded below.  Since $\Omega$ is
closed, $L_V(\im\Omega)=\dd(V\lrcorner \im\Omega)$.  Note that the
two differential forms $V\lrcorner \im\Omega$ and $\dd(V\lrcorner
\im\Omega)$ are defined on $X$ independently of $\alpha$.  From
this it follows that $\left|V\lrcorner \im\Omega|_{\M}\right|_0$
and $\left|\dd(V\lrcorner \im\Omega)|_{\M}\right|_0$ are bounded
independently of $\alpha$, where these are the induced norms on
$\M$. Thus
\begin{eqnarray*}
\left|\int_{\M}\left(\dd(V\lrcorner
\im\Omega|_{\M}\right)S\right|&\leq&
\left|\int_{\M}\left(\dd(V\lrcorner
\im\Omega|_{\M}\right)\varphi\bar{S}\right|+C\delta^{(n-2)/2}
\text{ by Lemma \ref{sbar}}\\
&=&\left|\int_{\M}(V\lrcorner
\im\Omega)\dd(\varphi\bar{S})\right|+C\delta^{(n-2)/2}\\
&\leq&\left|\int_{\M\cap B_{2\delta}}\left|V\lrcorner
\im\Omega|_{\M}\right|_0\!\cdot\!C|\dd\varphi|_0\right|+C\delta^{(n-2)/2}\\
&\leq&C\delta^{n-1}+C\delta^{(n-2)/2}.
\end{eqnarray*}
We now consider the second term.  Observe that since
$\im\dot{\Omega}$ is a form on $X$ defined independently of
$\alpha$, $|\im\dot{\Omega}|_{\M}|_0$ is bounded independently of
$\alpha$.  And obviously $|\im\dot{\Omega}|_{M_1\amalg
M_2}|_{0,M_1\amalg M_2}$ is bounded independently of $\alpha$.
Therefore
\begin{eqnarray*}
\int_{\M}(\im\dot{\Omega}|_{\M}) S &=&
\int_{\M}(\im\dot{\Omega}|_{\M})
\varphi\bar{S}+O(\delta^{(n-2)/2})\text{
by Lemma \ref{sbar}}\\
&=&\int_{M_1\amalg M_2}(\im\dot{\Omega}|_{M_1\amalg M_2})
\varphi\bar{S}
+O(\delta^{(n-2)/2})\\
&=&\int_{M_1\amalg M_2}(\im\dot{\Omega}|_{M_1\amalg M_2})
\bar{S}+O(\delta^{n/2})+O(\delta^{(n-2)/2})\text{ by Lemma
\ref{sbar}}\\
&=&\sqrt{{\vol(M_1)\vol(M_2)\over \vol(M_1)+\vol(M_2)}}
\left\langle {[M_1]\over\vol(M_1)}-{[M_2]\over \vol(M_2)},
[\im\dot{\Omega}]\right\rangle\\
& &+O(\delta^{(n-2)/2})\text{ by definition of }\bar{S}.
\end{eqnarray*}
Now the result follows from Lemma \ref{deform}.
\end{proof}

\section{The Full Linearized Deformation Operator}
We are now ready to prove the full injectivity estimate.
\begin{prop}[Full Injectivity Estimate]
For small enough $\nu$ independent of $\alpha$, $\DD
F_\alpha(0,0)$ satisfies the injectivity estimate
$C_I(\alpha)=C\epsilon^{-\nu}$ for sufficiently small $\alpha$.
That is,
$$
\|\DD F_\alpha(0,0)(u,a)\|_{\mathcal{B}'_\alpha}\geq
  {1\over C}\epsilon^{\nu}\|(u,a)\|_{\mathcal{B}_\alpha}.
$$
\end{prop}
\begin{proof}
\begin{eqnarray*}
|a|&\leq& C\left|\int_{\M} a\psi S\right| \text{ by the Key Lemma}\\
&=&C\left|\int_{\M}(\Delta u +a\psi)S\right|\text{ since }\Delta
u\text{ is
orthogonal to }S\\
&=&C\left|\int_{\M}\rho^2(\Delta u+a\psi)\rho^{-2}S\right|\\
&\leq&C\|\Delta u+a\psi\|_{\mathcal{B}'_\alpha}\int_{\M}|\rho^{-2}S|\\
&\leq&C\|\Delta u+a\psi\|_{\mathcal{B}'_\alpha}
\end{eqnarray*}
where the last line follows from the bounds on $|S|_0$ and
$\|\rho^{-1}\|_{L^2(\M)}$.  On the other hand,
\begin{eqnarray*}
\|u\|_{\mathcal{B}_{1,\alpha}}&\leq& C\epsilon^{-\nu}\|\Delta
u\|_{\mathcal{B}'_\alpha}\text{ by the Laplacian Injectivity Estimate}\\
&\leq&C\epsilon^{-\nu}(\|\Delta u+a\psi\|_{\mathcal{B}'_\alpha}+
\|a\psi\|_{\mathcal{B}'_\alpha})\\
&\leq&C\epsilon^{-\nu}\|\Delta u+a\psi\|_{\mathcal{B}'_\alpha}
\end{eqnarray*}
where the last line follows from the previous calculation and the
fact that $|\psi|_{C_\rho^{0,\beta}(\M)}$ is bounded independently
of $\alpha$, by the definition of $\psi$.  Finally, we deal with
the $Pu$ term.
\begin{eqnarray*}
\|(u,a)\|_{\mathcal{B}_\alpha}&\leq&C\epsilon^{-\nu}\|\Delta
u+a\psi\|_{\mathcal{B}'_\alpha}\text{ by combining the previous
two
calculations}\\
&\leq&C\epsilon^{-\nu}(\|\DD
F_\alpha(0,0)(u,a)\|_{\mathcal{B}'_\alpha}+\|Pu\|_{\mathcal{B}'_\alpha})\\
&\leq&C\epsilon^{-\nu}\|\DD
F_\alpha(0,0)(u,a)\|_{\mathcal{B}'_\alpha}+C\epsilon^{-\nu}\alpha^{1-\beta}
\|(u,a)\|_{\mathcal{B}_{\alpha}}\text{ by Lemma \ref{Pu}.}
\end{eqnarray*}
For sufficiently small $\nu$, $\epsilon^{-\nu}\alpha^{1-\beta}\to
0$ as $\alpha\to 0$.  So for small enough $\alpha$, we can absorb
the last term into the left-hand side.
\end{proof}

\begin{prop}
For small $\nu$, $\DD F_\alpha(0,0)$ is surjective for
sufficiently small $\alpha$.
\end{prop}
\begin{proof}
Consider the map $A:\mathcal{B}_\alpha\too\mathcal{B}'_\alpha$
defined by by $A:(u,a)\mapsto \Delta u + a\psi$.  By the proof of
the Full Injectivity Estimate together with the Key Lemma, it is
evident that $A$ is an isomorphism with $\|A^{-1}\|\leq
C\epsilon^{-\nu}$.  By Lemma \ref{Pu}, $\|P\|\leq
C\alpha^{1-\beta}$, therefore $\|A^{-1}P\|\leq C
\alpha^{1-\beta}\epsilon^{-\nu}$.  For small enough $\nu$,
$\alpha^{1-\beta}\epsilon^{-\nu}\to 0$ as $\alpha\to 0$, therefore
$I+A^{-1}P$ is invertible, and it follows that $A+P=\DD
F_\alpha(0,0)$ is surjective.
\end{proof}

\section{Solving the Deformation Problem}

The following Proposition can be found in \cite{ylee}.
\begin{prop}[Nonlinear Estimate]
For small $\alpha$, $F_\alpha$ satisfies a nonlinear estimate with
$C_N=C\epsilon^{-2}$ and $r_1= {1\over C}\epsilon^2$.  That is,
for $(h,t)\in\mathcal{B}_\alpha$ with
$\|(h,t)\|_{\mathcal{B}_\alpha}\leq r_1$,
$$\|\DD F_\alpha(h,t)(u,a)-\DD
F_\alpha(0,0)(u,a)\|_{\mathcal{B}'_\alpha}\leq C\epsilon^{-2}
\|(h,t)\|_{\mathcal{B}_\alpha}\!\cdot\!\|(u,a)\|_{\mathcal{B}_\alpha}$$
for all $(u,a)\in\mathcal{B}_\alpha$.
\end{prop}
The bound on $r_1$ is needed so that we can always assume that
$\tau=1$ in our definition of $F_\alpha$. Finally, we have the
following simple estimate.
\begin{prop}[Estimate of $F_\alpha(0,0)$]
For small enough $\alpha$,
$$\|F_\alpha(0,0)\|_{\mathcal{B}'_\alpha}\leq
C\alpha^3. $$
\end{prop}
\begin{proof}
Note that $F_\alpha(0,0)=\langle
\im\Omega,\vol_{\M}\rangle_{\M}=\sin\theta$.  Recall that
$\rho(x)\leq C \delta$ for $x\in \M\cap B_{\delta}$ and
$\delta={\alpha\over C_\delta}$. Then since $\sin\theta$ is
supported in $\M\cap B_{\delta}$,
\begin{equation}\label{Fest1}
|\rho^2\sin\theta|_0\leq C\delta^2|\sin\theta|_0\leq C\alpha^3
\end{equation}
where the second inequality follows from the bound on $\sin\theta$
from Lemma \ref{theta}.  Now we will estimate
$[\rho^{2+\beta}\sin\theta]_\beta$ by interpolation.  As in
(\ref{Fest1}), we see that
$$|\rho^{2+\beta}\sin\theta|_0\leq C\alpha^{3+\beta}.$$
We also have
\begin{eqnarray*}
|\nabla(\rho^{2+\beta}\sin\theta)|_0 &\leq&
|(2+\beta)\rho^{1+\beta}(\nabla\rho)\sin\theta|_0+|\rho^{2+\beta}
\nabla(\sin\theta)|_0\\
&\leq&C\delta^{1+\beta}|\sin\theta|_0
+C\delta^{2+\beta}|\nabla(\sin\theta)|_0\\
&\leq&C\alpha^{2+\beta}
\end{eqnarray*}
where the second line follows from the bound on $\nabla\rho$ and
the same reasoning as in (\ref{Fest1}), and the last line uses the
bounds on $\sin\theta$ and $\nabla(\sin\theta)$ from Lemma
\ref{theta}. Combining the two previous inequalities, we see that
$$[\rho^{2+\beta}\sin\theta]_\beta\leq C\alpha^3.$$
\end{proof}

Finally, let $r=(2C_I C_N)^{-1}={1\over C}\epsilon^{2+\nu}$, which
is less than $r_1$ for sufficiently small $\alpha$.  For small
enough $\nu$, we have ${r\over 2C_I}={1\over
C}\epsilon^{2+2\nu}>C\alpha^3\geq
\|F_\alpha(0,0)\|_{\mathcal{B}'_\alpha}$ for sufficiently small
$\alpha$.  We can now invoke the Inverse Function Theorem to find
a solution $F_\alpha(h,t)=0$ with
$\|(h,t)\|_{\mathcal{B}'_\alpha}\leq r$, and by elliptic
regularity, $h$ is smooth. Since $\|\nabla h\|_0\leq
C\epsilon^{1+\nu}$, it follows that there exists an embedded
special Lagrangian submanifold of $(X,J_t,\omega_t,\Omega_t)$ in a
$C\epsilon^{1+\nu}$-neighborhood of $\M$ for some $t<r$. Finally,
since the construction of $\M$ and $F_\alpha$ can be made to
depend smoothly on $\alpha$, and there is a unique solution to
$F_\alpha(h,t)=0$ in $B_r(0,0)$, we can also say that the embedded
SLags we constructed, as well as $t$, depend smoothly on $\alpha$.
This concludes the proof of the Main Theorem and the Torus Case.

\bibliographystyle{hamsplain}
\bibliography{research}

\end{document}